\documentclass[letterpaper, 10 pt, conference]{ieeeconf}
\IEEEoverridecommandlockouts

\usepackage{cite}
\usepackage{algorithm}
\usepackage{algpseudocode}
\usepackage{algorithmicx}
\usepackage{graphicx}
\usepackage{textcomp}

\usepackage{amssymb,amsmath,amsfonts,cases,mathtools}
\usepackage{microtype}
\usepackage{url}
\usepackage{hyperref}
\usepackage{xifthen}
\usepackage{xpatch}
\usepackage{amsthm}
\usepackage[dvipsnames]{xcolor}
\usepackage{tikz}
\usepackage{lipsum,booktabs}
\usepackage{pifont}
\usepackage{xcolor}

\newcommand{\n}{n}
\newcommand{\degree}{d}
\newcommand{\dii}{\degree_i}

\newcommand{\dd}{\mathbf{d}}

\newcommand{\cE}{\mathcal{E}}
\newcommand{\cG}{\mathcal{G}}

\newcommand{\cS}{\mathcal{S}}

\newcommand{\cN}{\mathcal{N}}

\newcommand{\cR}{\mathcal{R}}

\newcommand{\cV}{\mathcal{V}}

\newcommand{\cH}{\mathcal{H}}
\newcommand{\cC}{\mathcal{C}}
\newcommand{\cX}{\mathcal{X}}

\newcommand{\cCi}{\cC_i}
\newcommand{\cCj}{\cC_j}

\newcommand{\R}{\mathbb{R}}
\newcommand{\N}{\mathbb{N}}

\newcommand{\x}{x}
\newcommand{\s}{s}

\newcommand{\z}{z}
\newcommand{\y}{y}

\newcommand{\vv}{v}
\newcommand{\X}{X}

\newcommand{\pr}{\gamma}

\newcommand{\Ax}{A_{\x}}

\newcommand{\An}{A_{\nabla}}

\newcommand{\iter}{{t}}
\newcommand{\iterp}{{\iter + 1}}

\newcommand{\xt}{\x^{\iter}}

\newcommand{\zt}{\z^{\iter}}

\newcommand{\xtp}{\x^{\iterp}}

\newcommand{\ztp}{\z^{\iterp}}

\newcommand{\xstar}{\x^\star}

\newcommand{\xitp}{\x_{i}^{\iter+1}}

\newcommand{\zijtp}{\z_{ij}^{\iter+1}}

\newcommand{\xit}{\x_{i}^{\iter}}
\newcommand{\sit}{\s_{i}^{\iter}}

\newcommand{\zit}{\z_{i}^{\iter}}
\newcommand{\zijt}{\z_{ij}^{\iter}}
\newcommand{\zjit}{\z_{ji}^{\iter}}
\newcommand{\yit}{\y_{i}^{\iter}}

\newcommand{\xjt}{\x_{j}^{\iter}}
\newcommand{\sjt}{\s_{j}^{\iter}}

\newcommand{\yjt}{\y_{j}^{\iter}}

\newcommand{\bx}{\bar{\x}}
\newcommand{\bz}{\bar{\z}}

\newcommand{\px}{\x_\perp}
\newcommand{\pz}{\z_\perp}

\newcommand{\bxt}{\bx^{\iter}}
\newcommand{\bxtp}{\bx^{\iter+1}}
\newcommand{\bzt}{\bz^{\iter}}
\newcommand{\bztp}{\bz^{\iter+1}}

\newcommand{\pxt}{\px^{\iter}}
\newcommand{\pxtp}{\px^{\iter+1}}
\newcommand{\pzt}{\pz^{\iter}}
\newcommand{\pztp}{\pz^{\iter+1}}

\newcommand{\tpz}{\tilde{\z}_\perp}

\newcommand{\tpzt}{\tpz^\iter}
\newcommand{\tpztp}{\tpz^{\iterp}}

\newcommand{\Vpx}{W_{\px}}

\newcommand{\bb}{b}
\newcommand{\B}{B}
\newcommand{\M}{M}
\newcommand{\rr}{R}

\newcommand{\G}{G}

\newcommand{\one}{\mathbf{1}}

\newcommand{\oneNn}{\one_{N,\n}}

\newcommand{\zeq}{\z^{\text{eq}}}

\newcommand{\pzeq}{\pz^{\text{eq}}}

\newcommand{\norm}[1]{\left \|#1 \right \|}

\newcommand{\T}{^\top}
\newcommand{\inv}{^{-1}}

\newcommand{\blkdiag}{\text{blk\,diag}}
\newcommand{\until}[1]{\{1,\ldots,#1\}}
\newcommand{\col}{\textsc{col}}

\newcommand{\E}{\mathbb{E}}

\newcommand{\lip}{L}

\newcommand\oprocendsymbol{\hbox{$\square$}}
\newcommand\oprocend{\relax\ifmmode\else\unskip\hfill\fi\oprocendsymbol}
\def\eqoprocend{\tag*{$\square$}}
\def\er/{Erd\H{o}s-R\'enyi}

\newcommand{\m}{p}

\newcommand{\step}{\delta}

\newcommand{\msg}{m}

\newcommand{\msgijt}{\msg_{ij}^\iter}
\newcommand{\msgjit}{\msg_{ji}^\iter}

\newcommand{\msgijtp}{\msg_{ij}^{\iterp}}

\newcommand{\mt}{\msg^\iter}
\newcommand{\mtp}{\msg^{\iterp}}

\newcommand{\hm}{\hat{\msg}}

\newcommand{\hmjit}{\hm_{ji}^\iter}

\newcommand{\hmjitp}{\hm_{ji}^{\iterp}}

\newcommand{\hmt}{\hm^\iter}
\newcommand{\hmtp}{\hm^{\iterp}}

\newcommand{\cc}{c}

\newcommand{\RR}{\cR}

\newcommand{\tm}{\tilde{\msg}} 

\newcommand{\tmt}{\tm^\iter}
\newcommand{\tmtp}{\tm^{\iterp}}

\newcommand{\staz}{\X^\star}
\newcommand{\cstaz}{\cX^\star}

\newcommand{\mm}{m} 
\newcommand{\mij}{\mm_{ij}}
\newcommand{\mji}{\mm_{ji}} 

\newcommand{\drift}{g}

\newcommand{\Qu}{Q_{U}}
\newcommand{\Qw}{Q_{W}}

\newcommand{\lipc}{\lip_{\cC}}

\newcommand{\eigR}{r}

\newcommand{\xddots}{%
	\raise 4pt \hbox {.}
	\mkern 6mu
	\raise 1pt \hbox {.}
	\mkern 6mu
	\raise -2pt \hbox {.}
}

\graphicspath{{figs/}}

\newcommand{\as}{\text{ a. s.}}

\newcommand{\asequiv}{\!\stackrel{\text{a.s.}}{\equiv}\!}

\newcommand{\asleq}{\!\stackrel{\text{a.s.}}{\leq}\!}

\def\er/{Erd\H{o}s-R\'enyi}

\def\algo/{{ADMM-Tracking Gradient}}
\def\algob/{{ADMM-Tracking Gradient V2}}
\def\ralgo/{{Robust ADMM-Tracking Gradient}}
\def\ATG/{{ATG}}
\def\RATG/{{RATG}}
\def\GT/{{Gradient Tracking}}

\def\calgo/{{Compressed ADMM-Tracking Gradient}}

\def\citeDL/{\cite[Proposition~6.1.2]{bertsekas2015convex}}

\newcommand{\prsim}{0.1}
\newcommand{\stepsim}{0.5}
\newcommand{\rhosim}{0.9}
\newcommand{\alphasim}{0.9}

\newcommand{\fstar}{f^\star}

\newcommand{\captionspace}{-0em}
\newcommand{\blockspace}{-.65em}

\newcommand{\widfig}{.98}

\allowdisplaybreaks

\newtheorem{theorem}{Theorem}

\newtheorem{lemma}{Lemma}

\newtheorem{assumption}{Assumption}
\newtheorem{remark}{Remark}

\title{\LARGE \bf Modular Distributed Nonconvex Learning with Error Feedback} %

\author{Guido Carnevale, \IEEEmembership{Member, IEEE} Nicola Bastianello, \IEEEmembership{Member, IEEE} %
	\thanks{Guido Carnevale was partially supported by Fondi PNRR - Bando PE - Progetto PE11 - 3A-ITALY, ``Made in Italy Circolare e Sostenibile'' - Codice PE0000004, CUP: J33C22002950001.}
	\thanks{Nicola Bastianello was partially supported by the European Union's Horizon Research and Innovation Actions programme under grant agreement No. 101070162.}
	\thanks{Guido Carnevale is with the Department of Electrical,  Electronic and Information Engineering,  Alma Mater Studiorum - Universit\`a di Bologna,  Bologna, Italy. Email: {\tt\footnotesize{guido.carnevale@unibo.it}}}
	\thanks{Nicola Bastianello is with the School of Electrical Engineering and Computer Science, and Digital Futures, KTH Royal Institute of Technology, Stockholm, Sweden. Email: {\tt\footnotesize nicolba@kth.se}}}

\begin{document}

\maketitle
\thispagestyle{plain}
\pagestyle{plain}

\begin{abstract}
    In this paper, we design a novel distributed learning algorithm using stochastic compressed communications. 
	In detail, we pursue a modular approach, merging ADMM and a gradient-based approach, benefiting from the robustness of the former and the computational efficiency of the latter. 
	Additionally, we integrate a stochastic integral action (error feedback) enabling almost sure rejection of the compression error.
    We analyze the resulting method in nonconvex scenarios and guarantee almost sure asymptotic convergence to the set of stationary points of the problem.
	This result is obtained using system-theoretic tools based on stochastic timescale separation.
    We corroborate our findings with numerical simulations in nonconvex classification.
\end{abstract}


\section{Introduction}
The technological advances of the recent past have led to the adoption of devices with computation, sensing, and communication resources in a wide range of sectors, \textit{e.g.}, robotics, power grid, Internet-of-Things~\cite{nedic2018distributed,testa2023tutorial,molzahn2017survey,gao2021evaluation}.
These devices collect data from their environment, and, by relying on peer-to-peer communications, they can cooperatively learn from it, with the aim of training a more accurate model than any individual agent could.
In detail, the goal for $N$ agents is to solve the consensus (or cost-coupled) optimization problem
\begin{align}\label{eq:problem}
	\min_{\x \in \R^\n} f(x) \text{ with } f(x) := \sum_{i=1}^N f_i(\x),
\end{align}
where $f_i : \R^\n \to \R$ is a local, smooth, \textit{nonconvex} loss function, and $x \in \R^\n$ are the parameters of the model to be trained.
The trained models are often large in size ($n \gg 1$), requiring equally large communications; the solution is to utilize compressed communications \cite{richtarik_error_2024}.
This paper will focus on \textit{designing a distributed learning algorithm for nonconvex problems that employs compressed communications}.

The main classes of distributed optimization algorithms are based on \textit{gradient tracking}~\cite{nedic2017achieving} and Alternating Direction Method of Multipliers (\textit{ADMM})~\cite{makhdoumi2017convergence}.
Since gradient tracking schemes embed a linear average consensus block, they are not robust to asynchrony, packet losses, and generic additive disturbances (see~\cite{carnevale2025admm}).
On the other hand, ADMM suffers from a high computational cost since it requires the solution of an optimization problem at each iteration.
The so-called \algo/ method~\cite{carnevale2025admm} overcomes these limitations by inheriting computational efficiency from gradient-based optimization and resilience from consensus-ADMM~\cite{bastianello2022admm}.

The main weakness of \algo/ is its large communication burden.
Recent works in the literature combine \textit{compressed communications} and \emph{error feedback}~\cite{richtarik_error_2024} to tackle this issue in gradient tracking schemes~\cite{zhao2022beer,song_compressed_2022,xu_compressed_2025}.
In particular, \cite{song_compressed_2022} guarantees exact linear convergence for strongly convex problems, while~\cite{zhao2022beer,xu_compressed_2025} remove the convexity assumption.
These algorithms, however, employ standard average consensus to promote cooperation of the agents, which is not robust to non-idealities such as asynchrony and unreliable communications.
We remark that exact convergence with compressed communications can be achieved also by using dynamic scaling of the compression parameters~\cite{liao_robust_2024}.

Our main contribution is the development and the analysis of a novel distributed optimization algorithm named \calgo/.
The proposed method extends \algo/~\cite{carnevale2025admm} by embracing stochastic compression mechanisms, to reduce the communication burden, and error feedback, to preserve exact convergence.
The design of our method follows a modular approach where these new components are integrated into the consensus-ADMM framework inherited from \algo/.
By using system theory tools based on stochastic timescale separation, we establish almost sure asymptotic convergence to a set where the agents' estimates agree on a stationary point of the common nonconvex problem.
Indeed, unlike the strongly convex setting studied in~\cite{carnevale2025admm}, this paper considers nonconvex costs without the Polyak-Łojasiewicz condition. %
We also highlight that our proof strategy suitably exploits our modular design and stochastic timescale separation to conveniently isolate the compression and error feedback mechanisms' effects.
Thus, our approach paves the way for the design of a broad class of novel distributed schemes embedding these mechanisms.

\noindent\textit{Notation:}
the identity matrix in $\R^{m\times m}$ is $I_m$. 
The vectors of $N$ ones and zeroes are $1_N$ and $0_N$, while
$\oneNn := 1_N \otimes I_\n$ and $\mathbf{0}_{N,\n} := 0_N \otimes I_\n$ with $\otimes$ being the Kronecker product.
The cardinality of a finite set $S$ is $|S|$.
Both $\col(v_i)_{i \in \until{N}}$ and $\col(v_1,\dots,v_N)$ denote the vertical concatenation of the column vectors $v_1, \dots, v_N$.
We denote a graph as $\cG = (\mathcal{V}, \cE)$ where $\mathcal{V} = \until{N}$, $\cE \subseteq \mathcal{V} \times \mathcal{V}$, and $i$ and $j$ can exchange data only if $(i,j)\in\cE$. 
We use $\cN_i := \{j\in \mathcal{V} \mid (i,j) \in \cE\}$, $\dii:= |\cN_i|$, and $\dd := \sum_{i=1}^N \dii$. %
$R_+^n$ is the positive orthant of $R^n$.
The symbol $\blkdiag(M_i)_{i \in \until{N}}$ denotes the block-diagonal matrix whose $i$-th diagonal block is $M_i \in \R^{n_i \times n_i}$.

\section{Problem Description and Preliminaries}
\label{sec:setup}

\subsection{Problem Setup}

Our goal is to design a method to solve~\eqref{eq:problem} in a distributed way, namely, with update laws implementable over a network of $N$ agents using only local information and neighboring communication ruled by an undirect and connected graph $\cG = (\mathcal{V}, \cE)$.
We characterize problem~\eqref{eq:problem} as follows. %
\begin{assumption}\label{ass:cost}
    The function $f$ is differentiable and radially unbounded, while the gradients $\nabla f_i$ are $\lip$-Lipschitz continuous for some $\lip >0$ and all $i \in \cV$.\oprocend
\end{assumption}
Assumption~\ref{ass:cost} 
ensures the existence of $\fstar := \min_{x \in \R^{\n}}f(x)$ without imposing the uniqueness of stationary points of $f$.

\subsection{\algo/}
\label{sec:algo_design}

We now review \algo/~\cite{carnevale2023distributed,carnevale2025admm}, which serves as the foundation to the proposed approach.
Ideally, to solve~\eqref{eq:problem} each agent $i \in \cV$, at iteration $\iter \in \N$, could update its solution estimate $\xit \in \R^{\n}$ as 
\begin{align}
	\xitp = \xit + \pr\bigg(\frac{1}{N}\sum_{j=1}^N\xjt - \xit\bigg) - \pr\frac{\step}{N} \sum_{j=1}^N \nabla f_j(\xjt),\label{eq:desired_control_law}
\end{align}
where $\pr, \step > 0$ are tuning parameters.
In a distributed setup, the global quantities $\frac{1}{N}\sum_{j=1}^N \xjt$ and $\frac{1}{N}\sum_{j=1}^N \nabla f_j(\xjt)$ in~\eqref{eq:desired_control_law} are not available; thus \algo/ adopts consensus ADMM~\cite{bastianello2022admm} to robustly track them in a distributed way.
The key is to define the auxiliary optimization problem
\begin{align}\label{eq:cns_as_opt}
	\begin{split}
		\min_{
			\substack{
				(\y_1,\dots,\y_N) \in \R^{N\n}
				\\
				(\s_1,\dots,\s_N) \in \R^{N\n}
			}
		}
		&\frac{1}{2}\sum_{i=1}^N \left(\norm{\y_i - \xit}^2 + \norm{\s_i - \nabla f_i(\xit)}^2\right)
		\\
        &\text{s.t.:} \: \: \y_i = \y_j, \quad \s_i= \s_j  \: \forall (i,j)\in \cE,
	\end{split}
\end{align}
whose (unique) solution is $\frac{\mathbf{1}_{N,2\n}}{N}\sum_{j=1}^N \col(\xjt,\nabla f_j(\xjt))$, see~\cite{bastianello2022admm}.
Problem~\eqref{eq:cns_as_opt} matches the structure of the consensus problems addressed in~\cite{bastianello2020asynchronous} via the so-called Distributed ADMM method.
Notably, in~\eqref{eq:cns_as_opt}, the loss functions are quadratic and, thus, the ADMM subproblems admit a closed-form solution.
We then replace the unavailable $(\frac{1}{N}\sum_{j=1}^N \col(\xjt,\nabla f_j(\xjt)))$ in~\eqref{eq:desired_control_law} with the current solution estimates $(\yit,\sit) \in \R^{2\n}$ of problem~\eqref{eq:cns_as_opt} provided by the ADMM algorithm, giving rise to the overall distributed method reported in Algorithm~\ref{algo:algo}.
\begin{algorithm}[!ht]
	\begin{algorithmic}
		\State \textbf{Initialization}: $\x_i^0 \in \R^{\n}$, $\z_{ij}^0 \in \R^{2\n}$ for all $j \in \cN_i$. 
		\For{$\iter=0, 1, \dots$}
		\vspace{.1cm}
			\State $\begin{bmatrix}
				\yit
				\\
				\sit
			\end{bmatrix} = \frac{1}{1+\rho \dii}\left(\begin{bmatrix}\xit\\\nabla f_i(\xit)\end{bmatrix} + \sum_{j\in\cN_i} \zijt\right)$
			\vspace{.1cm}
			\State $\xitp = \xit + \pr(\yit - \xit) - \pr\step\sit$
			\For{$j \in \cN_i$}
				\State $\mij(\zijt,\yit,\sit) = -\zijt + 2\rho\begin{bmatrix}
					\yit
					\\
					\sit
				\end{bmatrix}$
				\State {\bf transmit} $\mij(\zijt,\yit,\sit)$ to $j$
				\State {\bf receive} $\mji(\zjit,\yjt,\sjt)$ from $j$
				\State $\zijtp = (1-\alpha)\zijt + \alpha\mji(\zjit,\yjt,\sjt)$
			\EndFor
		\EndFor
	\end{algorithmic}
	\caption{\algo/~\cite{carnevale2025admm}}%
	\label{algo:algo}
\end{algorithm}
We remark that the ADMM updates introduce the additional variables $\zijt \in \R^{2\n}$ maintained by agent $i$ for each of its neighbors $j \in \cN_i$, and that $\mij$ denotes the message sent from $i$ to $j \in \cN_i$.
Owing to this modular algorithm design, \algo/ demonstrates attractive robustness properties against imperfect communication and computation conditions, which it inherits from ADMM~\cite{carnevale2025admm}.
However, its main drawback is the high communication burden, especially when training a large model $\n \gg 1$.
To address this limitation, in the next section, we incorporate a communication compression mechanism.

\section{\calgo/}
\label{sec:calgo_design} 

In this section we extend Algorithm~\ref{algo:algo} to allow agents to use compressed communications and a local error feedback mechanism.
In detail, each agent $i \in \cV$, for each neighbor $j \in \cN_i$, implements a stochastic integral action by maintaining two variables $\msgijt, \hmjit \in \R^{2\n}$ and updating them as 
\begin{subequations}
    \begin{align}
        \label{eq:msgijt_update}
    	\msgijtp &= \msgijt + \cCi(-\zijt + 2\rho\col(\yit,\sit) - \msgijt),
        \\
        \label{eq:hmjit_update}
    	\hmjitp &= \hmjit + \cCj(-\zjit + 2\rho\col(\yjt,\sjt) - \msgjit),
    \end{align}
\end{subequations}
where $\cCi: \R^{2\n} \to \R^{2\n}$ generically denotes the randomized compression mechanism adopted by agent $i$.
Specifically, each $\msgijt$ implements a stochastic integral action asymptotically eliminating (with probability $1$) the error in the quantity $-\zjit + 2\rho\col(\yjt,\sjt)$ caused by compressed communication.
On the other hand, provided the initialization $\hm_{ji}^0 = \msg_{ji}^0$ is used, $\hmjit$ provides to agent $i$ a local copy of $\msgjit$ without requiring agent $j$ to transmit it via uncompressed communication.
The resulting distributed method is reported in Algorithm~\ref{algo:calgo}.
\begin{algorithm}[H]
	\begin{algorithmic}
		\State \textbf{Initialization}: $\x_i^0 \in \R^{\n}$, $\z_{ij}^0 \in \R^{2\n}$ and $\msg_{ij}^0 \in \R^{2\n}$ and $\hm_{ji}^0 = \msg_{ji}^0$ for all $j \in \cN_i$. 
		\For{$\iter=0, 1, \dots$}
		\vspace{.1cm}
			\State $\begin{bmatrix}
				\yit
				\\
				\sit
			\end{bmatrix} = \frac{1}{1+\rho \dii}\left(\begin{bmatrix}\xit\\\nabla f_i(\xit)\end{bmatrix} + \sum_{j\in\cN_i} \zijt\right)$
			\vspace{.1cm}
			\State $\xitp = \xit + \pr(\yit - \xit) - \pr\step\sit$
			\For{$j \in \cN_i$}
				\State {\bf transmit}  $\cCi(-\zijt + 2\rho\col(\yit,\sit) - \msgijt)$ to $j$
                \State {\bf receive} $\cCj(-\zjit + 2\rho\col(\yjt,\sjt) - \msgjit)$ from $j$
                \State 
				$\begin{aligned}
					\msgijtp &= \msgijt + \cCi(-\zijt + 2\rho\col(\yit,\sit) - \msgijt)
					\\
                 	\hmjitp &= \hmjit +\cCj(-\zjit + 2\rho\col(\yjt,\sjt) - \msgjit)
					\\
					\zijtp &= (1-\alpha)\zijt + \alpha\hmjit
				\end{aligned}$
			\EndFor
		\EndFor
	\end{algorithmic}
	\caption{\calgo/}%
	\label{algo:calgo}
\end{algorithm}
We underline that Algorithm~\ref{algo:calgo} uses only compressed communication.
Indeed, for each pair $(i,j) \in \cE$, once the compressed messages $\cCi(-\zijt + 2\rho\col(\yit,\sit) - \msgijt)$ and $\cCj(-\zjit + 2\rho\col(\yjt,\sjt) - \msgjit)$ are exchanged, the other algorithm updates are performed using only local variables.
The compression mechanisms $\cCi$ satisfy the following properties.
\begin{assumption}\label{ass:compression}
	For all $i \in \cV$, $\cCi$ is unbiased and contractive, namely, $\E[\cCi[v]] = v$ and $\E[\norm{\cCi(v) - v}^2] \leq (1 - \cc)\norm{v}^2$ for some $\cc \in (0,1)$ and all $v \! \in \! \R^{2\n}$. %
	Further, there exists $\lipc \!>\! 0$ such that $\cCi(v) \asleq \lipc\norm{v}$ for all $i \in \cV$ and $v \in \R^{\n}$.
	\oprocend
\end{assumption}
The next theorem proves almost sure convergence of the local estimates $\xt := \col(\xit)_{i \in \cV} \in \R^{N\n}$ generated by Algorithm~\ref{algo:calgo} to the set $\cstaz := \{x \in \R^{N\n} \mid x = \oneNn\xstar, \xstar \in \staz\}$, where $\staz := \{x \in \R^{n} \mid \nabla f(x) = 0\} \subset \R^{n}$.
\begin{theorem}\label{th:convergence}
	Let Assumptions~\ref{ass:cost} and~\ref{ass:compression} hold and consider the initialization in Algorithm~\ref{algo:calgo}.
	Then, there exist $\bar{\pr}, \bar{\alpha} > 0$ such that, for all $\rho >0$, $\alpha \in (0,\bar{\alpha})$, $\pr \in (0,\bar{\pr})$, $\step \in (0,\tfrac{2N}{\lip\sqrt{N}})$, the trajectories of Algorithm~\ref{algo:calgo} satisfy
	\begin{align*}
		\lim_{\iter\to\infty}\norm{\xt}_{\cstaz} = 0 \as.\eqoprocend
	\end{align*}
\end{theorem}
Theorem~\ref{th:convergence} is proved in Section~\ref{sec:proof_algo}.
The key intuition behind the proof is to interpret a suitable reformulation of Algorithm~\ref{algo:calgo} as a stochastic two-time-scale system (see~\cite{carnevale2024timescale}).
\begin{remark}%
	Algorithm~\ref{algo:algo} so far has been analyzed only in strongly convex scenarios~\cite{carnevale2025admm}. 
	However, Algorithm~\ref{algo:calgo} collapses to Algorithm~\ref{algo:algo} in the ``nominal'' case without compressed communication.
	Thus, Theorem~\ref{th:convergence} directly provides also the convergence proof of Algorithm~\ref{algo:algo} in nonconvex scenarios.
	\oprocend
\end{remark}

\section{Algorithm Analysis}
\label{sec:analysis_algo}

The analysis of \calgo/ follows three steps: 1) in Section~\ref{sec:system reformulation}, we interpret it as a stochastic two-time-scale system, i.e, the interconnection between a fast and slow stochastic subsystem; 2) in Sections~\ref{sec:bl_algo} and~\ref{sec:rs_algo}, we separately study the identified subsystems by relying on stochastic Lyapunov- and LaSalle-based arguments; 3) Section~\ref{sec:proof_algo} draws 1) and 2) together to prove Theorem~\ref{th:convergence}.

\subsection{Interpretation as Stochastic Two-Time-Scale System}
\label{sec:system reformulation}

\begin{figure}[H]
	\centering
	\vspace{\blockspace}
	\includegraphics[scale=1]{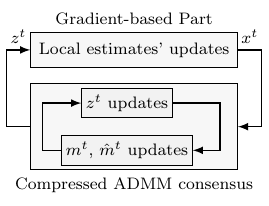}
	\vspace{\blockspace}
	\caption{Block diagram representing~\eqref{eq:algo_aggregate_form}.}
	\label{fig:block_diagram}
	\vspace{\blockspace}
\end{figure}

We start by providing the aggregate formulation of Algorithm~\ref{algo:calgo}.
To this end, we introduce the permutation matrix $P \in \R^{2\n\dd \times 2\n\dd}$ that swaps the $ij$-th and $ji$-th elements, and
\begin{align*}
	\Ax &:= \blkdiag\left\{ \col( \one_{\degree_i,\n}, \mathbf{0}_{\degree_i,\n} ) \right\}_{i \in \cV} \in \R^{2\n\dd \times N\n} \\
	\An &:= \blkdiag\left\{ \col( \mathbf{0}_{\degree_i,\n}, \one_{\degree_i,\n} ) \right\}_{i \in \cV} \in \R^{2\n\dd \times N\n} \\
    A &:= \blkdiag\left\{ \one_{\degree_i,2\n} \right\}_{i \in \cV} \in \R^{2\n\dd \times 2N\n} \\
    H &:= \blkdiag\left\{ \frac{1}{1+\rho\degree_i}I_{\n} \right\}_{i \in \cV}\in \R^{N\n \times N\n} \\
    \cH &:= \blkdiag\left\{ \frac{1}{1+\rho\degree_i}I_{2\n} \right\}_{i \in \cV} \in \R^{2N\n \times 2N\n}.
\end{align*}
By defining $\zit = \col(\zijt)_{j \in \cN_i} \in \R^{2\degree_i\n}$ and $\zt := \col(\zit)_{i \in \cV} \in \R^{2\n\dd}$ the aggregate formulation then is
\begin{subequations}\label{eq:algo_aggregate_form}
	\begin{align}
		\xtp &= \xt + \pr \left(H\left(\xt + \Ax\T\zt\right) -\xt\right)
		\notag\\
		&\hspace{.4cm}
		- \pr \step H  \left(\G(\xt) + \An\T\zt\right)\label{eq:algob_aggregate_form_x}
		\\
		\ztp &= (1-\alpha)\zt + \alpha P\hmt
		\label{eq:algob_aggregate_form_z}
		\\
		\mtp &= \mt + \cC(-\zt + 2\rho A\cH  (A\T\zt +\vv(\xt)) - \mt) 
		\label{eq:algob_aggregate_form_m}
        \\
        \hmtp &= \hmt + \cC(-\zt + 2\rho A\cH  (A\T\zt +\vv(\xt)) - \mt) 
		,\!\! \label{eq:algob_aggregate_form_mhat}
	\end{align}
\end{subequations}
in which $G(\x) \!:=\! \col(\nabla f_i(\x_i))_{i \in \cV}$, %
$\cC(z) \!:=\! \col(\cC_i(z_i))_{i \in \cV}$, %
and $\vv(\x) := \col(\col(\x_i,\nabla f_i(\x_i)))_{i \in \cV}$. %
Fig.~\ref{fig:block_diagram} reports a block diagram that graphically describes system~\eqref{eq:algo_aggregate_form}.
Under the initialization $\hm^0 = \msg^0$, we note that the trajectories of~\eqref{eq:algo_aggregate_form} satisfy $\hmt = \mt$ for all $\iter \in \N$. Therefore, from now on, we disregard~\eqref{eq:algob_aggregate_form_mhat} and focus on \eqref{eq:algob_aggregate_form_x},~\eqref{eq:algob_aggregate_form_z},~\eqref{eq:algob_aggregate_form_m}.
The next step is to rewrite these equations in a more convenient set of coordinates.
To this end, by using~\cite[Lemma~2]{bastianello2020asynchronous}, we observe that the matrix $P + 2\rho PA\cH  A\T$ has some eigenvalues in $1$ that are all semi-simple, while all the remaining ones lie in the open unit circle.
Then, let $\bb \in \N$ be the dimension of the subspace $\cS$ spanned by the eigenvectors of $P + 2\rho PA\cH  A\T$ associated to the eigenvalues equal to $1$, $\B \in \R^{2\n\dd \times \bb}$ be the matrix whose columns represent an orthonormal basis of $\cS$, and $\M \in \R^{2\n\dd \times \m}$ be such that $\M\T M = I$ and $\B\T \M = 0$, with $\m := 2\n\dd - \bb$.
The next lemma highlights useful properties of the matrices $\B$ and $\M$.
\begin{lemma}\label{lemma:results}
	[\cite[Lemma~IV.1]{carnevale2025admm}]
	The matrices $\B$ and $\M$ satisfy%
	\begin{subequations}\label{eq:results}
	\begin{align}
		\Ax\T \B = 0, \quad
		\An\T \B = 0, \B\T PA &= 0\label{eq:B_P_A}
		\\
		\B\T (I \! + \! P \! - \! 2\rho PA\cH A\T)&= 0.\label{eq:results_2}
	\end{align}
\end{subequations}
\end{lemma}
We also note that~\eqref{eq:B_P_A} and~\eqref{eq:results_2} imply
\begin{align}\label{eq:additional_results}
	\B\T P = -\B\T,
	\quad 
	\B\T P \B = - I. 
\end{align}
We then define the new coordinates $\bz := \B\T\z$, $\pz := \M\T\z$, and $\tm := \msg - (-I + 2\rho A\cH  A\T)\z - 2\rho A\cH\vv(\x)$, 
and, by using~\eqref{eq:results} and~\eqref{eq:additional_results}, we rewrite~\eqref{eq:algob_aggregate_form_x},~\eqref{eq:algob_aggregate_form_z},~\eqref{eq:algob_aggregate_form_m} as
\begin{subequations}\label{eq:algo_transformed}
	\begin{align}
		\xtp &= \xt + \pr \left(H\left(\xt + \Ax\T\M \pzt\right) -\xt\right)
		\notag\\
		&\hspace{.4cm}
		- \pr \step H  \left(\G(\xt) + \An\T\M \pzt\right)
		\label{eq:algo_transformed_px}
		\\
		\pztp &= (1-\alpha)\pzt + \alpha\M\T (- P + 2\rho PA\cH A\T)\M \pzt
		\notag\\
		&\hspace{.4cm}
		+ 2\alpha\rho \M\T PA\cH \vv(\xt) + \alpha M\T P\tmt
		\label{eq:algo_transformed_pz}
		\\
		\tmtp &= \tmt + \cC(- \tmt) 
		+\alpha\B\B\T P \tmt
		\notag\\
		&\hspace{.4cm}
		-2\rho PA\cH (\vv(\xtp) - \vv(\xt))
		\notag\\
		&\hspace{.4cm}
		-(-I + 2\rho A\cH  A\T)M(\pztp-\pzt),
		\label{eq:algo_transformed_m}
	\end{align}
\end{subequations}
and
$
	\bztp = \bzt + \alpha B\T P\tmt.
$
The obtained reformulation shows that $\bzt$ does not affect system~\eqref{eq:algo_transformed}, and, thus, we ignore it and focus on~\eqref{eq:algo_transformed} only.
In particular, we interpret~\eqref{eq:algo_transformed} as a stochastic two-time-scale system, where~\eqref{eq:algo_transformed_px} plays the role of \emph{slow} subsystem, while~\eqref{eq:algo_transformed_pz}-\eqref{eq:algo_transformed_m} is the \emph{fast} one.
Indeed, we can arbitrarily tune the variations of $\xt$ through $\pr$, while the fast state $\col(\pzt,\tmt)$ has equilibria $(\pzeq(\xt),0)$ parametrized in the slow state $\xt$, where $\pzeq: \R^{N\n} \to \R^{\m}$ is defined as 
\begin{align}
	\pzeq(x) &:=
	2\rho(I_{\m} - \RR)\inv \M\T PA\cH  \vv(x),
	\label{eq:equilibrium}
\end{align}
where we defined $\RR := -\M\T (P - 2\rho PA\cH A\T)\M \in \R^{\m \times \m}$.
Moreover, for all $x \in \R^{\n}$, it is possible to show that 
\begin{subequations}\label{eq:results_equilibrium}
	\begin{align}
		H \Ax\T \M \pzeq(x) &= {\oneNn\oneNn\T}x / {N} - H x
		\\
		H \An\T \M \pzeq(x) &= {\oneNn\oneNn\T} \G(x) / {N} -H \G(x).
	\end{align}
\end{subequations}
This interpretation allows us to proceed by separately analyzing the so-called boundary layer (originated from the fast part) and reduced (originated from the slow part) systems.

\subsection{Boundary Layer System}
\label{sec:bl_algo}

The boundary layer system is defined with an arbitrarily fixed $\xt \equiv x \in \R^{N\n}$ for all $\iter \in \N$ in the fast part~\eqref{eq:algo_transformed_pz}-\eqref{eq:algo_transformed_m} written in the error coordinate $\tpzt := \pzt - \pzeq(\x)$, namely%
\begin{subequations}\label{eq:bl_algo}
	\begin{align}
		\tpztp &= (1-\alpha)\tpzt + \alpha\RR\tpzt + \alpha M\T P\tmt
		\label{eq:bl_algo_pz}
		\\
		\tmtp &= \tmt + \cC(- \tmt) -\alpha(-I + 2\rho A\cH  A\T)P\tmt
		\notag\\
		&\hspace{.4cm}
		-\alpha(-I + 2\rho A\cH  A\T)M((-I + \RR)\tpzt).
		\label{eq:bl_algo_tm}
	\end{align}
\end{subequations}
The next lemma provides a function $U$ that ensures that the equilibrium solution $(\tpzt,\tmt) \asequiv 0$ to~\eqref{eq:bl_algo} is globally asymptotically stable in probability (see, e.g.,~\cite[Def.~2]{carnevale2024timescale}).
\begin{lemma}\label{lemma:bl_algo}
	There exist $\bar{\alpha} \in (0,1)$ and a radially unbounded and continuous function $U: \R^{(\m+2\dd\n)} \to \R_+$ such that, for all $\alpha \in (0,\bar{\alpha})$, the trajectories of~\eqref{eq:bl_algo} satisfy%
	\begin{subequations}\label{eq:U}
		\begin{align}
			&
			\!\!\E[U(\tpztp,\tmtp)] \! - \! U\!(\tpzt,\tmt)  \! \leq \! - b_1(\norm{\tpzt}^2 \!+\! \norm{\tmt}^2)\!\!\!
			\label{eq:U_1}
			\\
			&|U(\tpz,\tm) - U(\tpz^\prime,\tm^\prime)| 
			\leq b_2\norm{\begin{bmatrix}
				\tpz \! - \! \tpz^\prime
				\\\tm \! - \! \tm^\prime
			\end{bmatrix}}
			\left(
			\begin{bmatrix}
				\tpz 
				\\
				\tm
			\end{bmatrix} 
			\! - \!
			\begin{bmatrix}
				\tpz^\prime 
				\\
				\tm^\prime
			\end{bmatrix}
			\right)\!,
			\label{eq:U_2}
 		\end{align}
	\end{subequations}
	for all $\col(\tpzt,\tmt), \col(\tpz,\tm), \col(\tpz^\prime,\tm^\prime) \in \R^{(\m+2\dd\n)}$ and some $b_1, b_2 >0$.\oprocend
\end{lemma}
The proof of Lemma~\ref{lemma:bl_algo} is provided in Appendix~\ref{sec:lemma_bl}.

\subsection{Reduced System}
\label{sec:rs_algo}

Now, we study the reduced system obtained by considering $\pzt = \pzeq(\xt)$ into~\eqref{eq:algo_transformed_px} for all $\iter\in\N$, which by~\eqref{eq:results_equilibrium} yields
\begin{align}\label{eq:rs_algo}
	\hspace{-.1cm}
	\xtp \!&=\! \xt \!-\! \pr\big(I \!-\! \tfrac{1}{N}\oneNn\oneNn\T\big)\xt 
	\!-\! \tfrac{\pr\step}{N}\oneNn\oneNn\T\G(\xt).\!\!
\end{align}
The next lemma provides a Lyapunov function $W$ ensuring that the set $\cstaz$ is globally attractive for the trajectories of~\eqref{eq:rs_algo}.
To this end, let the matrix $\rr \in \R^{N\n \times (N-1)\n}$ be such that $\rr\T\oneNn = 0$ and $\rr\T\rr = I_{(N-1)\n}$, and $\phi: \R^{N\n} \to \R$ be
\begin{align}\label{eq:phi}
	\phi(x) := \norm{\begin{bmatrix}\oneNn\nabla f(\oneNn\T x/N)\T& (\rr x)\T\end{bmatrix}\T}.
\end{align}
We note that $\phi(x) = 0 \ \Leftrightarrow\  x \in \cstaz$.
\begin{lemma}\label{lemma:rs_algo}
	There exists a continuous and radially unbounded function $W: \R^{N\n} \to \R_+$ such that, for all $\pr \in (0,1)$ and $\step \in (0,2N/(\lip\sqrt{N}))$, the trajectories of~\eqref{eq:rs_algo} satisfy 
	\begin{subequations}\label{eq:W}
		\begin{align}
			W(\xtp) - W(\xt) &\leq -\step c_1\phi(\xt)^2
			\label{eq:W_1}
			\\
			|W(x + x^\prime) - W(x + x^{\prime\prime})| &\leq c_2\phi(x)\norm{x^\prime - x^{\prime\prime}} 
			\notag\\
			&\hspace{.4cm}
			+ c_2(\norm{x^\prime}^2 + \norm{x^{\prime\prime}}^2),
			\label{eq:W_2}
		\end{align}
	\end{subequations}
	for all $\xt, x, x^\prime, x^{\prime\prime} \in \R^{N\n}$ and some $c_1, c_2 > 0$.\oprocend
\end{lemma}
The proof of Lemma~\ref{lemma:rs_algo} is provided in Appendix~\ref{sec:lemma_rs}.

\subsection{Proof of Theorem~\ref{th:convergence}}
\label{sec:proof_algo}

To prove Theorem~\ref{th:convergence}, we apply the stochastic timescale separation result~\cite[Th.~1]{carnevale2024timescale}.
To this end, we need to verify:

(i) Lipschitz continuity of the dynamics and equilibrium function: verified by Lipschitz continuity of each $\nabla f_i$ (cf. Assumption~\ref{ass:cost}) and the bound on each $\cCi$ (cf. Assumption~\ref{ass:compression}).

(ii) Existence of a Lyapunov function proving global asymptotic stability in probability of the origin for the boundary layer system~\eqref{eq:bl_algo}: verified by Lemma~\ref{lemma:bl_algo}.

(iii) Existence of a Lyapunov function showing almost sure global attractiveness of $\cstaz$ for the reduced system~\eqref{eq:rs_algo} trajectories: verified by Lemma~\ref{lemma:rs_algo}.

We thus apply~\cite[Th.~1]{carnevale2024timescale} which ensures the existence of $\bar{\step} > 0$ such that, for all $\step \in (0,\bar{\step})$, the set $\{(x,\pz,\tm) \in \R^{N\n} \times \R^{\m} \times \R^{2\dd\n} \mid x\in\cstaz,\pz = \zeq(x), \tm = 0\}$ is almost surely globally attractive (cf.~\cite[Def.~1]{carnevale2024timescale}) for system~\eqref{eq:algo_transformed}.

\section{Numerical Simulations}
\label{sec:numerical_simualtions}

\begin{figure}[H]
	\centering
	\includegraphics[width=\widfig\linewidth]{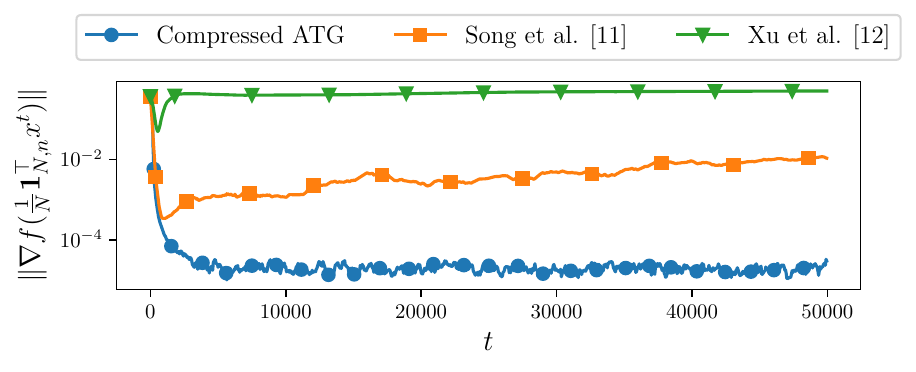}
	\vspace{\captionspace}
	\caption{Evolution of $\|\nabla f(\frac{1}{N}\oneNn\T\xt)\|$ achieved by Algorithm~\ref{algo:calgo},~\cite{song_compressed_2022}, and~\cite{xu_compressed_2025} with Top-1 compressors and in the presence of noise.}
	\label{fig:noise}
\end{figure}

In this section, we numerically test the effectiveness of \calgo/ on the nonconvex classification scenario addressed in~\cite{alghunaim2024local}. 
Specifically, in this scenario, each agent $i$ is equipped $m_i \in \N$ pairs of feature vectors and labels, and the corresponding problem reads as
\begin{align*}
	\min_{x \in \R^{\n}} \sum_{i=1}^N \frac{1}{m_i}\sum_{h=1}^{m_i}\log(1 + \exp(-b_{i,h}a_{i,h}\T x)) + \epsilon\sum_{\ell=1}^\n \frac{x_\ell^2}{1 + x_\ell^2},
\end{align*}
where $x_\ell$ is the $\ell$-th component of $x \in \R^\n$, and $a_{i,h} \in \R^\n$ and $b_{i,h} \in \{-1,1\}$ are the pairs of feature vector and label, while $\epsilon > 0$ is a regularization parameter.
We choose $N = 25$, $\n = 50$, $m_i = 250$, $\epsilon = 0.01$, a ring graph, while problem data are generated using \texttt{make\_classification} from \texttt{scikit-learn}~\cite{kramer2016scikit}.
We compare our method with those in~\cite{song_compressed_2022} and~\cite{xu_compressed_2025}. %
First, we test the case in which all the agents transmit a single component selected uniformly at random (Rand-1 compressor), namely, for all $i \in \cV$, we have
\begin{align}\label{eq:rand}
	\cCi(v) = v_{\bar{\ell}} e_{\bar{\ell}}, \quad \bar{\ell} \sim \mathcal{U}\until{\n},
\end{align}
where $\{e_1,\dots,e_\n\}$ is the canonical basis in $\R^\n$.
Then, we consider the case in which all the agents 
only send the largest component in absolute value (Top-1 compressor~\cite{islamov2023distributed}), i.e.,
\begin{align}\label{eq:top1}
	\cCi(v) = v_{\bar{\ell}} e_{\bar{\ell}}, \quad \bar{\ell} = \operatorname{arg\,max}_{\ell \in \until{\n}}|v_\ell|,
\end{align}
for all $i \in \cV$.
We empirically tune the parameters of the compared algorithms. For Alg.~\ref{algo:calgo} we set $\pr = \prsim$, $\step  = \stepsim$, $\rho = \rhosim$, and $\alpha = \alphasim$; for \cite{song_compressed_2022} we set $\step  = \stepsim$, $\beta = \gamma = \eta = 0.1$; and for \cite{xu_compressed_2025} we set $\step  = \stepsim$, $\gamma = \varphi_X = \varphi_Y = 0.1$.
\begin{figure}
	\centering 
	\includegraphics[width=\widfig\linewidth]{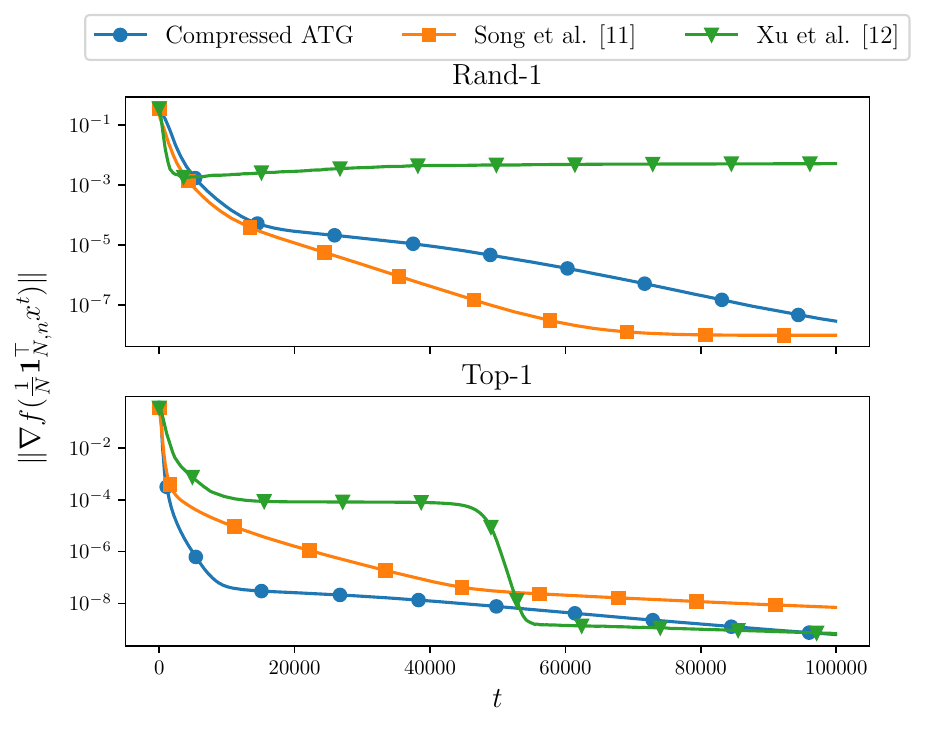}
	\vspace{\captionspace}
	\caption{Evolution of $\|\nabla f(\frac{1}{N}\oneNn\T\xt)\|$ achieved by Algorithm~\ref{algo:calgo}, \cite{song_compressed_2022}, and \cite{xu_compressed_2025}, with Rand-1 and Top-1 compressors.}
	\label{fig:comparison}
\end{figure}
Fig.~\ref{fig:comparison} reports the evolution of $\|\nabla f(\frac{1}{N}\oneNn\T\xt)\|$ for one run of the compared algorithms, using Rand-1 and Top-1 compression.
First of all, we remark that Alg.~\ref{algo:calgo} does indeed converge to a stationary point, as predicted by Theorem~\ref{th:convergence}.
Secondly, we notice that its performance, in terms of convergence rate, closely matches that of \cite{song_compressed_2022} in both cases. 
On the other hand, the algorithm in~\cite{xu_compressed_2025} converges exactly only by using Top-1 compressors.
We remark that since the considered problem is nonconvex, the algorithms' convergence rate is not uniform over $\iter$.
Table~\ref{tab:size-network} reports the approximate number of iterations needed by Alg.~\ref{algo:calgo} (with Top-1) to reach $\|\nabla f(\frac{1}{N}\oneNn\T\xt)\| \leq 10^{-6}$, for different $N$. We can see that larger networks have slower convergence speed, as they need to exchange more communications.
\begin{table}
	\vspace{-0em}
	\centering
	\caption{Num. of iterations to reach $\|\nabla f(\frac{1}{N}\oneNn\T\xt)\| \leq 10^{-6}$.}
	\vspace{-0em}
	\label{tab:size-network}
	\begin{tabular}{c|ccc}
	$N$ & $10$ & $25$ & $50$ \\
	\hline
	$t$ & $\sim 1600$ & $\sim 3300$ & $\sim 6100$ \\
	\end{tabular}
	\vspace{-0em}
\end{table}
We finally consider a noisy setting in which each agent $i$ sends to each neighbor $j \in \cN_i$ data in the form $d_{ij}^\iter + \eta^\iter_{ij}$, where $\eta^{\iter}_{ij}$ is drawn from a Gaussian distribution with zero-mean and variance $10^{-3}$, while $d_{ij}^\iter$ denotes the ``nominal'' data that would be sent from agent $i$ to $j$ at iteration $\iter$ in the noisy-free setting.
Fig.~\ref{fig:noise} shows the evolution of $\|\nabla f(\frac{1}{N}\oneNn\T\xt)\|$ in the considered schemes with Top-1 compressor and the same tuning parameters as before.
In detail, Fig.~\ref{fig:noise} shows that Alg.~\ref{algo:calgo} is more robust than the other two schemes.
Indeed, both methods in~\cite{song_compressed_2022} and~\cite{xu_compressed_2025} embed consensus schemes with a marginally stable subsystem that, differently from Alg.~\ref{algo:calgo}, affects the other system states.

\section{Conclusions}
\label{sec:conclusions}

We proposed a novel distributed method using compressed communications.
Following a modual approach, the proposed algorithm merges consensus-ADMM, gradient descent, and stochastic integral action. %
We relied on stochastic timescale separation to guarantee almost convergence to a stationary point in general 
scenarios without convexity assumptions.

\appendix

\subsection{Proof of Lemma~\ref{lemma:bl_algo}} 
\label{sec:lemma_bl}

Let us study the increment of the square norm $\norm{\tpz}^2$ along~\eqref{eq:bl_algo_pz}.
By adding $\pm\norm{(1-\alpha)\tpz + \alpha\RR\tpz}^2$, we obtain
\begin{align}
	&
	\norm{(1-\alpha)\tpz + \alpha\RR\tpz + \alpha \M\T P\tm}^2 - \norm{\tpz}^2
	\notag\\
	&=
	\norm{(1-\alpha)\tpz + \alpha\RR\tpz}^2 \! - \! \norm{\tpz}^2
	\! - \! 
	\norm{(1-\alpha)\tpz + \alpha\RR\tpz}^2
	\notag\\
	&\hspace{.4cm}
	+\norm{(1-\alpha)\tpz + \alpha\RR\tpz + \alpha \M\T P\tm}^2 
	\notag\\
	&\stackrel{(a)}{\leq} 
	(1-\alpha)\norm{\tpz}^2 + \alpha\norm{\RR\tpz}^2 - \alpha(1-\alpha)\norm{(\RR-I)\tpz}^2
	\notag\\
	&\hspace{.4cm}
	 - \norm{\tpz}^2
	+\alpha\norm{\M\T P\tm} \norm{(1-\alpha)\tpz + \alpha\RR\tpz}
	\notag\\
	&\hspace{.4cm}
	+ \alpha\norm{\M\T P\tm}\norm{(1-\alpha)\tpz + \alpha\RR\tpz + \alpha \M\T P\tm}
	\notag\\
	&\stackrel{(b)}{\leq} 
	- \alpha(1-\alpha)\eigR\norm{\tpz}^2 + \alpha2\norm{\tpz}\norm{\tm} + \alpha^2\norm{\tm}^2,
	\label{eq:increment_pz}
\end{align}
where $(a)$ uses standard properties of $\norm{\cdot}^2$, %
and $(b)$ uses the triangle inequality, $\norm{\M\T P\tm} \leq \norm{\tm}$, $\norm{(1-\alpha)\tpz + \alpha\RR}\leq\norm{\tpz}$, $\norm{\RR\tpz} \leq \norm{\tpz}$, and $-\norm{(\RR - I)\tpz}^2 \leq -\eigR\norm{\tpz}^2$ for some $\eigR > 0$, since $\RR$ has been derived by removing the part associated to unitary eigenvalues from an averaged operator~\cite{bastianello2020asynchronous}.
Now, let us define
\begin{align}
	&\drift(\tpz,\tm) 
	\! := \!
	(I - 2\rho A\cH  A\T)(M(-I + \RR)\tpz + P\tm),
	\label{eq:drift}
\end{align}
which allows us to rewrite~\eqref{eq:bl_algo_tm} as 
\begin{align}\label{eq:bl_algo_tm_compact}
	\tmtp = \tmt + \cC(- \tmt) + \alpha\drift(\tpzt,\tmt).
\end{align}
By using the triangle inequality, $\norm{(-I + 2\rho A\cH  A\T)M}\leq 1$, $\norm{P}\leq 1$, and $\norm{-I + \RR}\leq 2$, we get
\begin{align}\label{eq:lipschitz_drift}
	\norm{\drift(\tpz,\tm)} \leq \alpha2\norm{\tm} + \alpha2\norm{\tpz},
\end{align}
for all $(\tpz,\tm) \in \R^{\m} \times \R^{2\dd\n}$.
The expected increment of $\norm{\tm}^2$ along the trajectories of~\eqref{eq:bl_algo_tm_compact} reads as
\begin{align}
	&
	\E[\norm{\tm + \cC(- \tm) + \alpha\drift(\tpz,\tm)}^2] - \norm{\tm}^2
	\notag\\
	&\stackrel{(a)}{=} 
	\E[\norm{\tm + \cC(- \tm)}^2] - \norm{\tm}^2
	\notag\\
	&\hspace{.4cm}
	+\E[\norm{\tm + \cC(- \tm) + \alpha\drift(\tpz,\tm)}^2 -\norm{\tm + \cC(- \tm)}^2] 
	\notag\\
	&\stackrel{(b)}{\leq} 
	-\cc\norm{\tm}^2
	\notag
	+\alpha\E[\norm{\drift(\tpz,\tm)}\norm{\tm + \cC(- \tm) + \alpha\drift(\tpz,\tm)}]
	\notag\\
	&\hspace{.4cm}
	+\alpha\E[\norm{\drift(\tpz,\tm)}\norm{\tm + \cC(- \tm)}]
	\notag\\
	&\stackrel{(c)}{\leq} 
	-\cc\norm{\tm}^2 + \alpha4\norm{\tm}^2 + \alpha^2 4\norm{\tm}\norm{\tpz}
	\notag\\
	&\hspace{.4cm}
	+ \alpha4(1 + \alpha)\norm{\tm}^2 + \alpha^2 4\norm{\tpz}^2 + \alpha4(1 + \alpha)\norm{\tm}\norm{\tpz}
	\notag\\
	&= 
	-\cc\norm{\tm}^2
	+ \alpha4(1 + 2\alpha)\norm{\tm}^2 
	\notag\\
	&\hspace{.4cm}
	+ \alpha^2 4\norm{\tpz}^2 + \alpha4(1 + 2\alpha)\norm{\tm}\norm{\tpz},
	\label{eq:increment_tm}
\end{align}
where in $(a)$ we add $\pm\E[\norm{\tm + \cC(- \tm)}^2]$, in $(b)$ we use the fact that all the compressors are $\cc$-contractive to bound the first two terms and a standard property of the quadratic functions to bound the last ones,
while in $(c)$ we use the triangle inequality, the fact that the compressors are unbiased, and~\eqref{eq:lipschitz_drift}.
We now explicitly define the candidate Lyapunov function $U$ in~\eqref{eq:U} as 
$
	U(\tpz,\tm) := \norm{\tpz}^2 + \norm{\tm}^2. 
$
The chosen $U$ is quadratic and, thus, is continuous, radially unbounded, and verifies~\eqref{eq:U_2}.
In order to check also~\eqref{eq:U_1}, we use~\eqref{eq:increment_pz} and~\eqref{eq:increment_tm} to write 
\begin{align*}
	&
	\E[U((1 \! - \! \alpha)\tpz \! + \! \alpha\RR\tpz \! + \! \alpha \M\T\! P\tm,\tm \! + \! \cC(- \tm) \! + \! \drift(\tpz,\tm))]
	\notag\\
	&\hspace{.4cm}
	- U(\tpz,\tm)  
	\notag\\
	&\leq 
	-
	\begin{bmatrix}
		\norm{\tpz} 
		\\
		\norm{\tm}
	\end{bmatrix}\T 
	\underbrace{\begin{bmatrix}
		\alpha(\eigR - (4+\eigR)\alpha) & -\alpha(\alpha+2)
		\\
		-\alpha(\alpha + 2)& \cc - \alpha(4 + 9\alpha)
	\end{bmatrix}}_{\Qu(\alpha)}
	\begin{bmatrix}
		\norm{\tpz} 
		\\
		\norm{\tm}
	\end{bmatrix}.
\end{align*}
By Sylvester Criterion, there exists $\bar{\alpha} \in (0,1)$ such that, for all $\alpha \in (0,\bar{\alpha})$, it holds $\Qu(\alpha) > 0$ and the proof concludes.

\subsection{Proof of Lemma~\ref{lemma:rs_algo}} 
\label{sec:lemma_rs}

	We introduce a change of variables to isolate (i) the optimality and (ii) the consensus errors related to $\xt$.
	Namely, we introduce $\bx := \oneNn\T x/N \in \R^{\n}$ and $\px := \rr\T x\in \R^{(N-1)\n}$,
    and use them to rewrite~\eqref{eq:rs_algo} as the cascade system%
	\begin{subequations}\label{eq:rs_algo_explicit}
		\begin{align}
			\bxtp &= \bxt - \frac{\pr\step}{N}{\oneNn\T}\G(\oneNn\bxt + \rr\pxt)
			\label{eq:rs_algo_mu}
			\\
			\pxtp &= (1-\pr)\pxt,\label{eq:rs_algo_px}
		\end{align}
	\end{subequations}
    Now, let us define
    $
		\tilde{\G}\left(\bx,\px\right) := \frac{\oneNn\T}{N}(\G(\oneNn\bx)-\G(\oneNn\bx + \rr\px)),
    $
	and add $\pm\frac{\pr\step}{N} \nabla f(\bxt)$ in~\eqref{eq:rs_algo_mu}, thus obtaining
	\begin{subequations}\label{eq:sp_system_comapact}
		\begin{align}
			\bxtp &= \bxt - \frac{\pr\step}{N}\nabla f(\bxt) + \pr\step\tilde{\G}\left(\bxt,\pxt\right) 
			\label{eq:sp_system_comapact_mu}
			\\
			\pxtp &= (1-\pr)\pxt.\label{eq:sp_system_comapact_px}
		\end{align}
	\end{subequations}
	Once this formulation is available, we explicitly define $W$ as 
    \begin{align}\label{eq:W_rs_algo}
		W(\x) 
			&:=
		 f(\bx) - \fstar + \kappa/2\norm{\px}^2\!,\!
	\end{align}
	where $\kappa > 0$ will be fixed later.
	We note that $W(x)\ge0$ for all $x\in \R^{N\n}$ by definition of $\fstar$ and $W$ is continuous and radially unbounded since $f$ is radially unbounded and differentiable by Assumption~\ref{ass:cost}.
	To check~\eqref{eq:W_1}, we study the increment $\Delta f(\bxt,\pxt) \!:=\! f(\bxtp) \!-\! f(\bxt)$ along the trajectories of~\eqref{eq:sp_system_comapact_mu}.
	By Descent Lemma~\citeDL/, we get
    \begin{align}
    \Delta f(\bxt,\pxt) 
    &= 
	-\tfrac{\pr\step}{N}\norm{\nabla f(\bxt)}^2
    +\pr\step \nabla f(\bxt)\T\tilde{\G}(\bxt,\pxt)
    \notag\\
    &\hspace{.4cm}
    + \tfrac{\pr^2\step^2\lip\sqrt{N}}{N^2}( \nabla f(\bxt)\T \tilde{\G}(\bxt,\pxt)
    \notag\\
    &\hspace{.4cm}
	+  \tfrac{1}{2}\norm{\nabla f(\bxt)}^2 + \tfrac{N^2}{2}\|\tilde{\G}(\bxt,\pxt)\|^2).\!\!
		\notag\\
		&\stackrel{(a)}{\leq}
		-\tfrac{\pr\step}{N}\norm{\nabla f(\bxt)}^2
		+\pr\step \lip\sqrt{N}\norm{\nabla f(\bxt)}\norm{\pxt}
		\notag\\
		&\hspace{.4cm}
		\! + \!
		\tfrac{\pr^2\step^2\lip^2}{N^2}\!\norm{\nabla f(\bxt)}\!\norm{\pxt}
		+ \tfrac{\pr^2\step^2\lip^2}{2\sqrt{N}}\norm{\pxt}^2
		\notag\\
		&\hspace{.4cm}
		+\tfrac{\pr^2\step^2\lip\sqrt{N}}{2N^2} \norm{\nabla f(\bxt)}^2,
		\label{eq:DeltaV_mu_second}
	\end{align}
	where $(a)$ uses 
    $
		\|\tilde{\G}(\bxt,\pxt)\| \!\leq\! \lip\norm{\pxt}/\sqrt{N}
    $, which follows by Lipschitz continuity of $\nabla f_i$ and $\norm{\rr} \!=\! 1$.
	Now, by evaluating the increment $\Delta \Vpx(\pxt) := \frac{\kappa}{2}\norm{\pxtp}^2 - \frac{\kappa}{2}\norm{\pxt}^2$ along the trajectories of system~\eqref{eq:sp_system_comapact_px}, we get
	\begin{align}
		\Delta \Vpx(\pxt) &= -\kappa\pr(1 - \pr/2)\norm{\pxt}^2.
		\label{eq:DeltaV_px}
	\end{align} 
	Since $\nabla f$ is Lipschitz continuous, $\nabla W$ is Lipschitz continuous too.
	Further, it holds $\nabla W(x) = \frac{1}{N}\oneNn\nabla f(\oneNn\T x/N) + \kappa \rr\T x$ $\Rightarrow$ $\norm{\nabla W(x)} \leq \max\{1/N,\kappa\}\sqrt{2}\phi(x)$ (see the definition of $\phi$ in~\eqref{eq:phi}).
	Thus, condition~\eqref{eq:W_2} can be verified by Descent Lemma~\citeDL/.
	As for condition~\eqref{eq:W_1}, we define $\Delta W(\xt) := \Delta f(\bxt,\pxt) + \Delta \Vpx(\pxt)$ and, by using~\eqref{eq:DeltaV_mu_second} and~\eqref{eq:DeltaV_px}, the trajectories of~\eqref{eq:sp_system_comapact} satisfy
	\begin{align}
		\Delta W(\xt) 
		&\leq - \pr \begin{bmatrix}
			\norm{\nabla f(\bxt)}
			\\
			\norm{\pxt}
		\end{bmatrix}\T\! \Qw(\pr,\step)\!
		\begin{bmatrix}
			\norm{\nabla f(\bxt)}
			\\
			\norm{\pxt}
		\end{bmatrix}\!
		,\label{eq:V_muV_px}	
	\end{align}
	where %
	\begin{align*}
		\Qw(\pr,\step) 
		\!\!:=\!\! 
		\begin{bmatrix}
			\step/N -\frac{\pr\step^2\lip\sqrt{N}}{2N^2}& -\frac{\step}{2}\!(\lip\sqrt{N} \! + \! \frac{\pr\step\lip\sqrt{N}}{2N^2}L)
			\\
			-\frac{\step}{2}\!(\lip\sqrt{N} \! + \! \frac{\pr\step\lip\sqrt{N}}{2N^2}L)& \kappa(1 - \pr/2) - \frac{\pr\step^2\lip^2}{2\sqrt{N}}
		\end{bmatrix}
		\!\!.
	\end{align*}
	Note that the upper-left block of $\Qw(\pr,\step)$ and $\kappa(1 - \pr/2)$ (in the lower-right block) are positive for all $\step \in (0,2N/(\lip\sqrt{N}))$ and $\pr \in (0,1)$, while the cross entries of $\Qw(\pr,\step)$ do not depend on $\kappa$.
	Thus, there exists $\bar{\kappa} > 0$ such that, for all $\kappa > \bar{\kappa}$, we have $\Qw(\pr,\step) > 0$ and the proof follows.

\bibliographystyle{IEEEtran} 
\bibliography{biblio}

\end{document}